\documentclass[12pt,psamsfonts]{article}
\usepackage[dvips]{epsfig}
\usepackage{graphics}
\usepackage{amsmath}
\usepackage{amsfonts}
\usepackage{amssymb}
\usepackage{amsthm}
\usepackage{bm}
\usepackage{enumerate}
\usepackage[mathscr]{eucal}
\usepackage{inputenc}
\usepackage[T2A]{fontenc}
\theoremstyle{plain}
\newtheorem{theorem}{Theorem}[section]
\newtheorem{lemma}{Lemma}[section]
\newtheorem{definition}{Definition}[section]
\newtheorem{remark}{Remark}[section]

\newtheorem*{nonumtheorem}{Theorem}
 \numberwithin{equation}{section}

\begin{document}
\title{\textbf{The truncated Fourier operator. V.
}}
\author{\textbf{Victor Katsnelson, Ronny Machluf}}
\date{\ }%

 \maketitle
 \thanks{\hspace*{-2.0ex}\textbf{Mathematics Subject
Classification: (2000).} 35S30, 43A90.\endgraf
\hspace*{-2.0ex}\textbf{Keywords:}
\begin{minipage}[t]{0.70\linewidth}{\small Truncated Fourier operator,
spectral analysis, prolate spherical functions.}\end{minipage}

 \abstract{The operator Fourier truncated on a finite symmetric interval is considered.
 The limiting behavior of its spectrum is discussed as the length of the interval
 tends to infinity.}
\setcounter{section}{4}
\section{\!\!\!\!Spectral theory of the Fourier operator\\ truncated on a
finite symmetric interval.}
In this section we discuss the behavior of eigenvalues \(\lambda_k(a)\) of the
truncated Fourier operator \(\mathscr{F}_E\) as functions of \(a\):
\begin{equation}
\mathscr{F}_{[-a,a]}e_k=\lambda_ke_k\,.
\end{equation}
\hspace*{1.0ex}\textbf{1.}\,
As we have seen in Section 4, (see \cite{KaMa2}), the operator \(\mathscr{F}_{[-a,a]}\)
commutes with the differential operator \(\mathcal{L}\) which is
generated by the differential expression
\begin{equation}
Lx(t)=-\frac{d\,}{dt}\Bigg(\bigg(1-\frac{t^2}{a^2}\bigg)\frac{dx}{dt}\Bigg)+t^2x(t)
\end{equation}
on the set of functions \(x\)
which are defined on the interval \((-a,a)\) and satisfied the boundary
conditions
\begin{equation}
\lim_{t\to{}-a+0}(t+a)\frac{dx(t)}{dt}=0,\,\quad
\lim_{t\to{}a-0}(t-a)\frac{dx(t)}{dt}=0\,.
\end{equation}
This operator \(\mathcal{L}\) appeared as the operator \(\mathcal{L}_I\) in the
section 4, and is the member of the family \(\lbrace\mathcal{L}_U\rbrace\)
which consists of all selfadjoint differential operators generated by the formal differential
operator \(L\). Here \(U\) is an arbitrary \(2\times2\) unitary matrix indexing the
family. See \cite[Definition 4.19 and Theorem 4.1]{KaMa2}.
 The operator \(\mathcal{L}_I\) is the member of the family corresponding
to the matrix \(U=I\), where \(I\) is the identity matrix.
Since we do not deal with the operator \(\mathcal{L}_U\) other than the operator
\(\mathcal{L}_I\), we omit the subindex \(I\) and use the notation
\(\mathcal{L}\) instead of \(\mathcal{L}_I\).

The operator \(\mathcal{L}\) is a selfadjoint operator in
\(L^2([-a,a])\) with discrete spectrum. The spectrum   of the
operator \(\mathcal{L}\) consists of the sequence of eigenvalues
\(\mu_k(a)\) which tend to \(+\infty\). Each of these eigenvalues
is of multiplicity one. (These well known facts facts are
formulated as \cite[Theorem 4.3]{KaMa2}.)

 We enumerate the eigenvalues of \(\mathcal{L}\)
in the increasing order:
\begin{equation}%
\label{EVpPa}
0=\mu_0(a)<\mu_1(a)<\mu_2(a)<\,\,\ldots\,\,<\mu_k(a)<\,\,\ldots\,\,.
\end{equation}
 The appropriate eigenfunction \(e_k(t,a)\) can be normalized such that it
depends on \(a\) smoothly. Indeed the eigenvalue problem for the operator \(\mathcal{L}\)
is of the form
\begin{gather}
-\frac{d\,}{dt}\Bigg(\bigg(1-\frac{t^2}{a^2}\bigg)\frac{de_k(t,a)}{dt}\Bigg)+t^2e_k(t,a)=
\mu_k(a)e_k(t,a),\quad -a<t<a\,,\notag\\
|e_k(-a,a)|<\infty,\quad |e_k(a,a)|<\infty\,.
\label{EVPA}
\end{gather}
Changing variable \(t\to{}at\), we come to the eigenvalue problem
\begin{subequations}
\label{EVP}
\begin{gather}
-\frac{d\,}{dt}\Bigg(\big(1-t^2\big)\frac{dg_k(t,a)}{dt}\Bigg)+a^2t^2g_k(t,a)=
\gamma_k(a)g_k(t,a),\quad-1<t<1\,,
\label{EVP1a}\\[0.2ex]
|g_k(-1,a)|<\infty,\quad |g_k(1,a)|<\infty\,,
\label{EVP1b}
\end{gather}
\end{subequations}
where for \(a\geq0\)
\begin{equation}%
\label{EVpP}
0=\gamma_0(a)<\gamma_1(a)<\gamma_2(a)<\,\,\ldots\,\,<\gamma_k(a)<\,\,\ldots\,\,,
\end{equation}
and for \(a>0\)
\begin{equation}%
\label{MwN}
\gamma_k(a)=a^2\mu_k(a)\,.
\end{equation}%
The  eigenvalue  problems of the family \eqref{EVPA} are considered
on the \emph{different} intervals \((-a,a)\). In other words, the appropriate
operators act in the spaces \(L^2(-a,a)\) which depend on \(a\).

The eigenvalue problems \eqref{EVP} are considered in the single interval
\((-1,1)\). The dependence on \(a\) is transferred to the operator: before the
potential \(t^2\) the factor \(a^2\) appears.

Since the eigenvalues of the boundary value problem \eqref{EVP} are simple,
and the differential operator \eqref{EVP1a} depends on the parameter \(a\)
analytically,
 each eigenvalue \(\gamma_k(a),\,k=0,1,2,\,\,\ldots\,\,,\) depends on
\(a\) analytically for real \(a\). This is a standard result of the perturbation theory.

\begin{definition}
\label{PSpF}
The eigenfunctions \(g_k(t,a)\) of the boundary value problem
\eqref{EVP} are said to be\,\footnote%
{Various authors disagree not only on notation for these functions,
but also in their method of normalization.
}
 the \emph{angular prolate spherical functions.}
\end{definition}

 Since the equation \eqref{EVP1a} is invariant with respect to
the change of variable \(t\to-t\) and the eigenvalue \(\nu_k(a)\)
simple, the eigenfunction \(g_k(t,a)\) is either even, or odd with
respect to \(t\). It turns out, that\,%
\footnote{Recall that the eigenvalues and eigenfunctions of the
boundary value problem \eqref{EVP} are enumerated according to
\eqref{EVpP}. } %
 for even \(k\) the eigenfunction
\(g_k(t,a)\) is even, and for odd \(k\) the eigenfunction
\(g_k(t,a)\) is odd with respect to \(t\). This can be proved
studying the oscillation properties of solutions of the
differential equation \eqref{EVP1a}.

Different normalization of the prolate spherical functions are used in literature.
In the monograph \cite{Fl} the normalization
\begin{subequations}
\label{NoBC}
\begin{equation}
\label{NoBCe}
g_k(t,a)_{|_t=0}=P_k(t)_{|_t=0},
\quad \frac{dg_k(t,a)}{dt}_{|_{t=0}}=0\,,
\quad k\textup{ \ is even},
\end{equation}
\begin{equation}
\label{NoBCo}
g_k(t,a)_{|_{t=0}}=0,\quad
 \frac{dg_k(t,a)}{dt}_{|_{t=0}}=\frac{dP_k(t)}{dt}_{|_{t=0}}\,,
 \quad k\textup{ \ is odd}\,.
\end{equation}
\end{subequations}
where \(P_k(t)\) are the Legendre polynomials,
\begin{equation}
\label{LegP}
P_k(t)=\frac{1}{2^{k}k!}\frac{d^k\,}{dt^k}(t^2-1)^k\,.
\end{equation}
In particular,
\begin{subequations}
\label{SpV}
\begin{alignat}{3}
\label{SpVe}
P_k(t)_{|_{t=0}}&=&  \ &\frac{(-1)^{\frac{k}{2}}k!}{2^k\big(\frac{k}{2}\big)\!!\,\big(\frac{k}{2}\big)\!!}\qquad ,
\quad &k\text{ is even},\\
\label{SpVeo}
\frac{dP_k(t)}{dt}_{|_{t=0}}&=&\ \  \ &
\frac{(-1)^{\frac{k-1}{2}}(k+1)!}{2^k\big(\frac{k-1}{2}\big)\!!\,\big(\frac{k+1}{2}\big)\!!},
\quad &k\text{ is odd}\,.
\end{alignat}
\end{subequations}

Under this normalization, each of the eigenfunctions \(g_k(t,a)\) is the
solution of the Cauchy problem for the differential equation \eqref{EVP1a}
with the initial conditions \eqref{NoBC} posed in the point \(t=0\).
The boundary conditions \eqref{EVP1b} are satisfied automatically
since the number \(\gamma_k(a)\), which appears in the right hand side of \eqref{EVP1a},
is an eigenvalue. The function \(g_k(t,a)\) depends on \(a\) analytically
for real \(a\)
because this function is a solution of the differential equation \eqref{EVP1a} whose
coefficients \(a^2t^2\) and \(\gamma_k(a)\) depends on \(a\) analytically
and the initial conditions \eqref{NoBC} do not depend on \(a\).

For \(a=0\), the differential equation \eqref{EVP1a} is the Legendre
differential equation,  eigenfunctions \(g_{k}(t,0)\)
of the eigenvalue problem \eqref{EVP}
are the Legendre polynomials, \eqref{LegP},
\begin{equation}
\label{LegPo}
g_k(t,0)=P_k(t),
\end{equation}
 and the eigenvalues \(\gamma_k(0)\) are the numbers
\begin{equation}
\label{LegPEv}
\gamma_k(0)=k(k+1)\,.
\end{equation}
Since the set of functions
\(\lbrace{}g_k(t,a)\rbrace_{k=0,\,1,\,2,\,\,\ldots\,\,}\) is the
set of all eigenfunctions of a selfadjoint differential operator
with discrete spectrum, this set forms an orthogonal basis in the
Hilbert space \(L^2([-1,1])\). The eigenfunctions \(e_k(t,a)\) of
the boundary value problem \eqref{EVPA}, in other words the
eigenfunctions of the operator \(\mathcal{L}\), are
\begin{equation}
\label{ReEF}%
 e_k(t,a)=g_{k}(ta^{-1},a)\,.
\end{equation}
The set of functions
\(\lbrace{}e_k(t,a)\rbrace_{k=0,\,1,\,2,\,\,\ldots\,\,}\) forms an
orthogonal basis in the Hilbert space \(L^2([-a,a])\).

\noindent
 \hspace*{1.0ex}\textbf{2.}\,The selfadjoint differential operator
\(\mathcal{L}\) commutes with the truncated Fourier operator
\(\mathscr{F}_{[-a,a]}\) as well as with the adjoint operator
\(\mathscr{F}_{[-a,a]}^{\,\,\ast}\).
 (This fact is
formulated in \cite{KaMa2}) as Theorem 4.2.) Thus the operator
\(\mathcal{L}\) commutes also with the operator
\(\mathscr{F}_{[-a,a]}^{\,\,\ast}\mathscr{F}_{[-a,a]}\).

 Since the eigenvalues
of the operator \(\mathcal{L}\) are simple, the eigenfunctions
\(e_k(t,a)\) of the operator \(\mathcal{L}\) also are an
eigenfunctions of each of the operators \(\mathscr{F}_{[-a,a]}\),
\(\mathscr{F}_{[-a,a]}^{\,\,\ast}\) and
\(\mathscr{F}_{[-a,a]}^{\,\,\ast}\mathscr{F}_{[-a,a]}\):
\begin{subequations}
\label{EFTrF}%
\begin{equation}
\label{EFTrFo}%
 \mathscr{F}_{[-a,a]}e_k(\,.\,,a)=\lambda_k(a)e_k(\,.\,,a)\,,\quad
 (k=0,\,1,\,2,\,\ldots\,\,)\,.
\end{equation}
\begin{equation}
\label{EFTrFa}%
 \mathscr{F}_{[-a,a]}^{\,\,\ast}e_k(\,.\,,a)=\overline{\lambda_k(a)}e_k(\,.\,,a)\,,\quad
 (k=0,\,1,\,2,\,\ldots\,\,)\,.
\end{equation}
\begin{equation}
\label{EFTrFao}%
 \mathscr{F}_{[-a,a]}^{\,\,\ast}\mathscr{F}_{[-a,a]}e_k(\,.\,,a)
 =|\lambda_k(a)|^2e_k(\,.\,,a)\,,\quad
 (k=0,\,1,\,2,\,\ldots\,\,)\,.
\end{equation}
\end{subequations}
Here \(\lambda_k(a)\) is the eigenvalue of the operator
\(\mathscr{F}_{[-a,a]}\) corresponding to the eigenfunction
\(e_k(\,.\,,a)\). Since the operator \(\mathscr{F}_{[-a,a]}\) is
normal, the numbers \(\overline{\lambda_k(a)}\) and
\(|\lambda_k(a)|^2\) are the eigenvalues of the operators
\(\mathscr{F}_{[-a,a]}^{\,\,\ast}\) and
\(\mathscr{F}_{[-a,a]}^{\,\,\ast}\mathscr{F}_{[-a,a]}\)
corresponding to the eigenfunction \(e_k(\,.\,,a)\).

 In this section, we discuss the distribution of the
eigenvalues \(\lambda_k(a)\) and their behavior, mainly for large
\(a\). We already know that the eigenvalues \(\lambda_k(a)\) are
located within the cross which is the union of the horizontal
interval \([-1,1]\) and the vertical interval \([-i,i]\). The
eigenvalues \(\lambda_k(a)\) corresponding to even (w.r.t. \(t\))
eigenfunctions \(e_k(t,a)\) belong to \([-1,1]\), the eigenvalues
\(\lambda_k(a)\) corresponding to odd (w.r.t. \(t\))
eigenfunctions \(e_k(t,a)\) belong to \([-i,i]\). (The even
functions \(e_k(t,a)\) are eigenfunctions of the cosine-transform,
the odd functions \(e_k(t,a)\) are eigenfunctions of the
sine-transform. See formula (2.14) in \cite{KaMa1}.) We already
mentioned that the functions \(e_k(t,a)\) with even and odd
indices \(k\) are even and odd functions of \(t\) respectively.
The operator
\(\mathscr{F}_{[-a,a]}^{\,\,\ast}\mathscr{F}_{[-a,a]}\) is an
integral operator:
\begin{equation}%
\label{FStF}%
(\mathscr{F}_{[-a,a]}^{\,\,\ast}\mathscr{F}_{[-a,a]}x)(t)=
\frac{1}{\pi}\int\limits_{-a}^{a}\frac{\sin{}a(t-\tau)}{t-\tau}\,x(\tau)\,d\tau\,.
\end{equation}%
The eigenvalue problem for the operator
\((\mathscr{F}_{[-a,a]}^{\,\,\ast}\mathscr{F}_{[-a,a]}x)\) is
\begin{equation}
\label{OEVP}
\frac{1}{\pi}\int\limits_{-a}^{a}\frac{\sin{}a(t-\tau)}{t-\tau}\,x(\tau)\,d\tau=%
\sigma(a)x(t)\,,\quad -a\leq{}t\leq a\,.
\end{equation}
Changing variable \(at\to{}t\) and \(t\to\dfrac{t}{a}\) we see that
the the eigenvalues \(\sigma(a)\) of the problem \eqref{OEVP} are the same that the eigenvalues
 of each of the problems
\begin{equation}
\label{OEVP1}
\frac{1}{\pi}\int\limits_{-a^2}^{a^2}\frac{\sin{}(t-\tau)}{t-\tau}\,x(\tau)\,d\tau=%
\sigma(a)x(t)\,,\quad -a^2\leq{}t\leq a^2\,,
\end{equation}
and
\begin{equation}
\label{OEVP2}
\frac{1}{\pi}\int\limits_{-1}^{1}\frac{\sin{a^2}(t-\tau)}{t-\tau}\,x(\tau)\,d\tau=%
\sigma(a)x(t)\,,\quad -1\leq{}t\leq 1\,.
\end{equation}
In the problem \eqref{OEVP1}, the kernel
\(\dfrac{1}{\pi}\dfrac{\sin{}(t-\tau)}{t-\tau}\) of the integral
operator does not depend on \(a\), but the interval \([-a^2,a^2]\)
depends on \(a\). In the  problem \eqref{OEVP2}, the interval
\([-1,1]\) does not depend on \(a\) but the kernel
\(\dfrac{1}{\pi}\dfrac{\sin{a^2}(t-\tau)}{t-\tau}\) depends on
\(a\). In  \cite[Section VI, pp.\,59-61]{SlPo} it is proved that
the eigenvalues of the problem \eqref{OEVP2} are pairwise
distinct. A transparent presentation of this result can be found
also in \cite{KPS}, Chapter 2, \S 2.2, pp.\,55\,-\,56. (In
\cite{KPS}, the eigenvalue problem \eqref{OEVP2} is considered).
In the proof of the non-degeneracy property of eigenvalues of the
problem \eqref{OEVP}-\eqref{OEVP1}-\eqref{OEVP2} the fact that the
eigenfunctions of the problem \eqref{OEVP} for the integral
operator are the same that the eigenfunctions \(e_k(t,a)\) of the
problem \eqref{EVPA} for the differential operator is essentially
used.

The integral operator
\begin{equation*}
x(t)\to\frac{1}{\pi}\int\limits_{-1}^{1}\frac{\sin{a^2}(t-\tau)}{t-\tau}\,x(\tau)\,d\tau
\end{equation*}
is of the form
\begin{equation*}
\frac{1}{\pi}\int\limits_{-1}^{1}\frac{\sin{a^2}(t-\tau)}{t-\tau}\,x(\tau)\,d\tau=
P_1\mathscr{F}^{\ast}P_{a^2}\mathscr{F}P_{{1}_{\big|_{L^2([-1,1])}}}\,,
\end{equation*}
where \(\mathscr{F}\) is the (non-truncated) Fourier operator acting in
\(L^2((-\infty,\infty)\), and for positive \(c\), \(P_c\) is the orthogonal projector
 from
\(L^2((-\infty,\infty))\) onto \(L^2((-c,c)\)\,.
Since for \(0<c_1<c_2\), \(P_1\mathscr{F}^{\ast}P_{c_1}\mathscr{F}P_{1}\leq{}P_1\mathscr{F}^{\ast}P_{c_2}\mathscr{F}P_{1}\),
then from minimax principle of Courant it follows that the eigenvalues \(\sigma_k(a)\)
of the problem \eqref{OEVP2} behave monotonically:
\begin{equation*}
\sigma_k(a_1)\leq{}\sigma_k(a_2)\quad\textup{if}\ \ a_1<a_2\,.
\end{equation*}
In this stage of the reasoning, the numbers \(\sigma_k(a)\) are enumerated
according to the rule
\begin{multline*}
\frac{1}{\pi}\int\limits_{-1}^{1}\frac{\sin{a^2}(t-\tau)}{t-\tau}g_k(t,a)=
\sigma_{k}(a)g_k(t,a), \ \ -1\leq{}t\leq{}1\,,\\[-2.0ex]
g_k(t,a) \ \ \textup{are the eigenfunctions of the problem \eqref{EVP}}\,.
\end{multline*}
After the non-degeneracy and the monotonicity properties of the eigenvalues
\(\sigma_{k}(a)\) of the operator \(\mathscr{F}_{[-a,a]}^{\,\,\ast}\mathscr{F}_{[-a,a]}\)
are established, the ordering
\begin{equation}
\label{OrdSi}
\sigma_{0}(a)>\sigma_{1}(a)>\sigma_{2}(a)>\,\,\ldots\,\,>\sigma_{k}(a)>\,\,\ldots\,,\quad
\textup{for every}\ \ a>0
\end{equation}
can be obtained from an analysis of the asymptotic
behavior of \(\sigma_{k}(a)\) for small \(a\). Such analysis was first done in \cite[Section VI]{SlPo}.
In this analysis, the property \eqref{LegPo} is essentially used.
In \cite{Wid} it is shown that for fixed \(k\),
\begin{equation}
\label{AsFkSp}
\sigma_{k}(a)=
2\pi\Big(\frac{a^2}{4}\Big)^{2k+1}(k!\,)^{-2}\big(1+o(1)\big),\ \ \textup{as}\ \ a\to+0.
\end{equation}
Since \(\sigma_k(a)=\big|\lambda_k(a)\big|^2\), the ordering
\begin{equation}
\label{OrdLa}
|\lambda_{0}(a)|>|\lambda_{1}(a)|>|\lambda_{2}(a)|>\,\,\ldots\,\,>|\lambda_{k}(a)|>\,\,\ldots\,,\quad
\textup{for every}\ \ a>0
\end{equation}
of the absolute
values of the eigenvalues \(\lambda_k(a)\) of the truncated Fourier transform
\(\mathscr{F}_{[-a,a]}\) holds as well.
The values of the arguments of the complex numbers \(\lambda{k}(a)\) also cam be
obtained  analysing the asymptotic behavior of \(\lambda{k}(a)\) for small \(a\).
For fixed \(k\),
\begin{equation}
\label{AsLsa}
\lambda_{k}(a)=i^{k}\sqrt{2\pi}\bigg(\frac{a}{2}\bigg)^{2k+1}\big(k!\,\big)^{-1}
\big(1+o(1)\big),\ \ \textup{as}\ \ a\to+0.
\end{equation}
In particular, for fixed \(k\) and \(a\to+0\), the argument \(\arg\lambda_k(a)\) tends
to one of the numbers \(\arg1,\,\arg{}i,\,\arg(-1),\,\arg{}(-i)\) depending
on the residue of \(k\) by \(\mod4\). Since the argument of the eigenvalues \(\lambda{k}(a)\)
can can take only the values \(\arg1,\,\arg{}i,\,\arg(-1),\,\arg{}(-i)\), then
\begin{equation}
\label{Argl}
\lambda_k(a)=i^{k}\big|\lambda_k(a)\big|\,.
\end{equation}
for \(a\) positive and small enough.
The eigenvalue \(\lambda_k(a)\) depends on \(a\) continuously. By increasing of \(a\),
the absolute value \(\big|\lambda_k(a)\big|\) increases. In particular, as \(a\)  increases,
the eigenvalue \(\lambda_k(a)\) is separated from the point \(\lambda=0\).
Therefore the eigenvalue \(\lambda_k(a)\) can not jump from the interval
\([0,i^k]\) to another interval \([0,i^m]\), where \(m\not=k\,\mod4\).
Thus, \emph{the equality \eqref{Argl} hods for every \(a>0\)}, and not only for small \(a\).
In \cite{Fu} is shown that for fixed \(k\), the asymptotic equality
\begin{equation}
\label{AsKfai}
1-\sigma_{k}(a)=4\sqrt{\pi}8^k\big(k!\big)^{-1}a^{2k+1}e^{-2a^2},\ \ \textup{as}\ \ a\to\infty\,.
\end{equation}
Therefore,
\begin{equation}
\label{AsKfail}
1-\big|\lambda_{k}(a)\big|=2\sqrt{\pi}8^k\big(k!\big)^{-1}a^{2k+1}e^{-2a^2},\ \ \textup{as}\ \ a\to\infty\,.
\end{equation}
We summarize the above stated facts as
\begin{theorem}
\label{BEFk}{\ }\\
\hspace*{1.5ex}\textup{1.}\,
\begin{minipage}[t]{0.95\linewidth}
 The eigenfunctions \(e_k(t,a)\) of the truncated Fourier operator \(\mathscr{F}_{[-a,a]}\)
 are expressed in  terms of the angular prolate spherical functions \(g_k(t,a)\)
 \textup{(Definition \ref{PSpF})}:
 \begin{equation*}%
 e_k(t,a)=g_k(ta^{-1}\!,\,a)\,.
 \end{equation*}
\end{minipage}\\[2.0ex]
\hspace*{1.5ex}\textup{2.}\,
\begin{minipage}[t]{0.95\linewidth}
 For every \(a\), the absolute values \(|\lambda_k(a)|\) of
the eigenvalues  \(\lambda_k(a)\) of the truncated Fourier operator \(\mathscr{F}_{[-a,a]}\)
are pairwise different, and their ordering \eqref{OrdSi}-\eqref{EFTrFo} agrees with the
ordering \eqref{EVpPa}-\eqref{EVPA} of the eigenvalues \(\mu_k(a)\) of the
differential operator \(\mathcal{L}\).
\end{minipage}\\[2.0ex]
\hspace*{1.5ex}\textup{3.}\,
\begin{minipage}[t]{0.95\linewidth}
The argument of the eigenvalues  \(\lambda_k(a)\) is:
\begin{equation*}%
\lambda_k(a)=i^k|\lambda_k(a)|\,.
\end{equation*}
\end{minipage}\\[2.0ex]
\hspace*{1.5ex}\textup{4.}\,
\begin{minipage}[t]{0.95\linewidth}
If \(k\) is fixed,
and \(a\) increases from
\(0\) to \(+\infty\), then the eigenvalue \(\lambda_k(a)\) moves
monotonically from the point
\(\lambda=+0\) to the point \(\lambda=i^k\).
\end{minipage}\\[2.0ex]
\end{theorem}
Theorem \ref{BEFk} is related to the behavior of the individual
eigenvalue \(\lambda_k(a)\) as a function of \(a\).
The next result is related to the behavior of the set
\(\lbrace\lambda_k(a)\rbrace_{k=0,\,1,\,2,\,\ldots}\) of all eigenvalues
for large \(a\). The next Theorem claims that for large \(a\), the set
\(\lbrace\lambda_k(a)\rbrace_{k=0,\,1,\,2,\,\ldots}\) fills, in a sense, the whole cross
with the endpoints \(1,\,i,\,-1,\,-i\).

If \(K(t,\tau),\,-a\leq{}t,\tau\leq{}a\,,\) is a smooth kernel, and \(\textbf{K}\) is
an integral operator in \(L^2([-a,a])\), then the operators \(\textbf{K}\)  is a trace class operator, and
\begin{equation*}
\textup{trace}\,\textbf{K}=\int\limits_{-a}^{a}K(\xi,\xi)\,d\xi\,,
\end{equation*}
The operator \(\textbf{K}^2\) also is a trace class operator, and its trace also may
expressed as
\begin{equation*}
\textup{trace}\,(\textbf{K}^2)=\int\limits_{-a}^{a}\int\limits_{-a}^{a}K(\xi,\eta)K(\eta,\xi)\,d\xi\,d\eta\,.
\end{equation*}
If \(\textbf{K}\) is the operator: \(\textbf{K}=\mathscr{F}_{[-a,a]}^{\,\,\ast}\mathscr{F}_{[-a,a]}\),
that is the integral operator with the kernel \[K(t,\tau)=\dfrac{1}{\pi}\dfrac{\sin{}a(t-\tau)}{t-\tau},\]
then \(K(\xi,\xi)=a\), and the trace can be calculated explicitly:
\begin{subequations}
\label{tr}
\begin{equation}
\label{trK}
\textup{trace}\,\mathscr{F}_{[-a,a]}^{\,\,\ast}\mathscr{F}_{[-a,a]}=\frac{2}{\pi}\,a^2\,.
\end{equation}
Trace of the operator \(\big(\mathscr{F}_{[-a,a]}^{\,\,\ast}\mathscr{F}_{[-a,a]}\big)^2\)
can not be calculated explicitly, but can be estimated from below:
\begin{equation}
\label{trKs}
\textup{trace}\,\big(\mathscr{F}_{[-a,a]}^{\,\,\ast}\mathscr{F}_{[-a,a]}\big)^2\geq
\frac{2}{\pi}\,a^2-\frac{2}{\pi^2}\ln^+a-1\,.
\end{equation}
\end{subequations}
The estimate \eqref{trKs} was obtained in \cite[formula (4.2)]{LaP2}.

For a trace class operator, its trace is equal to the sum of its eigenvalues.
(This fact is attributed to V.B.\,Lidskii.)
Thus, if \(\sigma_k(a)\) are eigenvalues of the operator
\(\mathscr{F}_{[-a,a]}^{\,\,\ast}\mathscr{F}_{[-a,a]}\), which are ordered in
the decreasing order, \eqref{OrdSi}, then
\begin{subequations}
\label{EVI}
\begin{equation}
\label{EVI1}
\sum_k{\sigma_k(a)}=\frac{2}{\pi}a^2
\end{equation}
and
\begin{equation}
\label{EVI2}
\sum_k{(\sigma_k(a))^2}\geq\frac{2}{\pi}a^2-\frac{2}{\pi^2}\ln^+a-1\,.
\end{equation}
Moreover, since the operator \(\mathscr{F}_{[-a,a]}^{\,\,\ast}\mathscr{F}_{[-a,a]}\)
is contractive,
\begin{equation}
\label{EVI3}
0<\sigma_k(a)<1\,.
\end{equation}
\end{subequations}
Subtracting the inequality \eqref{EVI2} from the equality \eqref{EVI1}, we obtain
the inequality
\begin{equation}
\label{ITZ}
\sum\limits_{k=0}^{\infty}\sigma_k(a)(1-\sigma_k(a))\leq\frac{2}{\pi^2}\ln^+a+1\,
\end{equation}
In view of \eqref{EVI3}, every summand in the left hand side of \eqref{ITZ} is positive.

Let a number \(\varepsilon\) is given, \(0<\varepsilon<1/2\). From \eqref{ITZ} it follows that
\begin{equation*}
\varepsilon\,\cdot\,\sum\limits_{k:\,\varepsilon<\sigma_k}(1-\sigma_k(a))
\leq\frac{2}{\pi^2}\ln^+a+1\,,
\end{equation*}
and
\begin{equation*}
\sum\limits_{k:\,\varepsilon<\sigma_k}1\leq\sum\limits_{k:\,\varepsilon<\sigma_k}\sigma_k(a)
+\varepsilon^{-1}\Big(\frac{2}{\pi^2}\ln^+a+1\Big)\,.
\end{equation*}
Taking into account \eqref{EVI1}, we come to the inequality
\begin{subequations}
\label{Reg}
\begin{equation}
\label{ZeReg}
\#\lbrace{}k:\,\sigma_k(a)>\varepsilon\rbrace<\frac{2}{\pi}a^2+
\varepsilon^{-1}\Big(\frac{2}{\pi^2}\ln^+a+1\Big)\,.
\end{equation}
In the same way, from \eqref{ITZ} it follows that
\begin{equation*}
\varepsilon\,\cdot\!\!\!\!\sum\limits_{k:\,\sigma_k\leq{}1-\varepsilon}\sigma_k(a)
\leq\frac{2}{\pi^2}\ln^+a+1\,,
\end{equation*}
Taking into account the equality \eqref{EVI1}, we come to the inequality
\begin{equation*}
\sum\limits_{k:\,\sigma_k>1-\varepsilon}\sigma_k(a)
\geq\frac{2}{\pi}a^2-\varepsilon^{-1}\Big(\frac{2}{\pi^2}\ln^+a+1\Big)\,.
\end{equation*}
In view of \eqref{EVI3}, the inequality
\begin{equation}
\label{OneReg}
\#\lbrace{}k:\,\sigma_k(a)>1-\varepsilon\rbrace\geq\frac{2}{\pi}a^2-
\varepsilon^{-1}\Big(\frac{2}{\pi^2}\ln^+a+1\Big)\,.
\end{equation}
holds.
Since \(\lambda(1-\lambda)\geq{}\varepsilon(1-\varepsilon)\) for
\(\lbrace\lambda:\,\varepsilon\leq{}\lambda\leq{}1-\varepsilon{}\rbrace\),
it follows from \eqref{ITZ} that
\begin{equation}
\label{IntReg}
\#\lbrace{}k:\,\varepsilon\leq\sigma_{k}(a)\leq{}1-\varepsilon\rbrace\leq
\frac{1}{\varepsilon(1-\varepsilon)}\Big(\frac{2}{\pi^2}\ln^+a+1\Big)\,.
\end{equation}%
\end{subequations}

We summarize the inequalities \eqref{Reg}:\\[1.0ex]
\hspace*{0.03\linewidth}%
\begin{minipage}[t]{0.92\linewidth}
\emph{For large \(a\), the eigenvalues \(\sigma_k(a)\)
of the operator \(\mathscr{F}_{[-a,a]}^{\,\,\ast}\mathscr{F}_{[-a,a]}\),
arranged in the decreasing order, first are very
close to one, then are very close to zero. The transition from the values
equal almost one   to
the values equal almost zero occurring  in the interval of values \(k\) which is centered
at the point \(k=\frac{2}{\pi}a^2\) and grows in width at the rate of only
\(\log{}a\).}
\end{minipage}\\[2.0ex]

 This result, even if it is not formulated in the explicit form,
as well as  its derivation, which is based on the estimates \eqref{tr} of the
 traces of the operators \(\mathscr{F}_{[-a,a]}^{\,\,\ast}\mathscr{F}_{[-a,a]}\)
 and \((\mathscr{F}_{[-a,a]}^{\,\,\ast}\mathscr{F}_{[-a,a]})^2\), appear first
 in \cite{LaP2}.

However, this result is not enough for our goal. We need a more detailed information
on the distribution of the eigenvalues \(\sigma_k(a)\) for large \(a\) and the indices
\(k\) belonging to the "transition" interval
\begin{equation}
\label{TrInt}
k:\,\,\,\frac{2}{\pi}a^2-\varepsilon^{-1}\Big(\frac{2}{\pi^2}\ln^+a+1\Big)\leq{}k\leq
\frac{2}{\pi}a^2+\varepsilon^{-1}\Big(\frac{2}{\pi^2}\ln^+a+1\Big)\,.
\end{equation}
The appropriate result was formulated
(as a conjecture) in the paper \cite[page 106]{Sl1} of D.\,Slepian. The conjecture of
D.\,Slepian is:

 Let \(\delta\) is the root of smallest absolute value of the equation
 \begin{equation*}
\frac{2}{\pi}a^2+\frac{2}{\pi}\log(2a)-\arg\Gamma\bigg(\frac{1}{2}+i\frac{\delta}{2}\bigg)=
k-\frac{1}{2}\,,
 \end{equation*}
 where \(\Gamma\) is the Euler gamma function and \(\arg\Gamma\big(\frac{1}{2}\big)=0\).
 Then for large \(k\) and \(a\), the approximate equality
\begin{equation*}
\sigma_k(a)\approx(1+e^{\pi\delta})^{-1}
\end{equation*}
holds.

 In particular, if \(b\) is fixed and
 \begin{subequations}
\label{IntAs}
\begin{equation}
\label{IntAs1}
k(a,b)=\Big[\frac{2}{\pi}(a^2+b\ln{}(2a))\Big]\,,
\end{equation}
where the brackets \([\,\,]\) denotes "largest integer in", then
\begin{equation}
\label{IntAs2}
\lim_{a\to\infty}\sigma_{k(a,b)}=(1+e^{\pi{}b})^{-1}\,.
\end{equation}
\end{subequations}
Though no rigorous mathematical proof of this approximation result is done
in \cite{Sl1}, the approximation is justified by numerical examples there.

The first rigorous proof of the asymptotic relation \eqref{IntAs} was done
in \cite{ClMe}. In particular, see formulas (3.27),\,(3.28) there.
The proof of \eqref{IntAs} given in \cite{ClMe} is based on the thorough study
of the asymptotic behavior of the prolate spherical functions \(g_k(t,a)\) in the
interval \(-1\leq{}t\leq{}1\) for large value of the parameter \(a\) and \emph{all}
\(k\). The method of the study can be classified as a hard analysis method.
The method used involves asymptotic solving the differential equation which defines
the prolate functions by using the classical WKB approximation.
Later H.\,Landau and H.\,Widom, \cite{LaWi}, propose the proof based on the theory of
integral operators with difference kernel \(k(t-\tau)\), where \(k\)
is a rapidly decreasing function.
In \cite[Theorem 1]{LaWi}, the eigenvalue distribution of the appropriate integral operator is
derived from the Fourier transform of the function \(k\), in accordance with
the classical Szego method. In \cite[Theorem 2]{LaWi}, even more general result is obtained
related to the integral operator with the kernel \(\frac{1}{\pi}\frac{\sin{}a^2(t-\tau)}{t-\tau}\)
considered on the set which is a finite union of non-intersecting intervals.

The analysis of the proof of the asymptotic result \eqref{IntAs} shows
that the limiting relation \eqref{IntAs} holds not only for fixed \(b\), but
for \(b\) from any fixed finite interval of the real axis. Moreover, the limit
in \eqref{IntAs2} is uniform with respect to \(b\) from any fixed finite interval:

\begin{nonumtheorem}[Cloizeaux-Mehta,\,\cite{ClMe};\,\,Landau-Widom,\,\cite{LaWi}]
Given positive numbers \(N\) and \(\varepsilon\), there exists the number %
\(A=A(\varepsilon,\,N)\) such that for any \(a\) satisfying the inequality
\(a\geq{}A\) and for any \(b\) belonging to the interval
\([-N,\,N]\), the inequality
\begin{equation}
\label{appr}
\big|\sigma_{k(a,b)}-(1+e^{\pi{}b})^{-1}\big|\leq\varepsilon\,,
\end{equation}
where the index \(k(a,b)\) is defined in \eqref{IntAs1}\,.
\end{nonumtheorem}
This result is formulated explicitly neither in the paper \cite{ClMe}, nor in the paper
\cite{LaWi}.
However the proof of this result can be extracted from any of this papers.

Let \(\kappa\), \(\kappa<1\), be an arbitrary positive number. In fact we  assume that \(\kappa\)
is  small.
We choose a positive number \(N=N(\kappa)\) such that
\begin{equation}%
\label{ChooN}
\big(1+e^N\big)^{-1}<\kappa/2\,.
\end{equation}%
Let \(I\) be an arbitrary interval of length \(\kappa\) which is contained
in the interval \([0,1]\): \(|I|=\kappa,\,I\subset[0,1]\,.\)
\emph{We fix this interval.}
 Since
the function \(b\to(1+e^b\big)^{-1}\) is monotonic for \(b\in(-\infty,\infty)\), and
\((1+e^b\big)^{-1}_{\,\,|_{b=\infty}}=0,\ \ \ %
(1+e^b\big)^{-1}_{\,\,|_{b=-\infty}}=1\), there exists the unique value of \(b\)
such that the point \((1+e^b\big)^{-1}\) is the center point of the interval \(I\).
\emph{Let us fix this \(b\) and denote it by \(b(I)\).}
Since the smallest possible value and the largest possible value for
the center point of the interval \(I\) of length \(\kappa\) are the points
\(\kappa/2\) and \(1-\kappa/2\), the number \(b(I)\) satisfy the inequality
\begin{equation}%
\label{IncN}%
-N\leq{}b(I)\leq{}N\,,
\end{equation}%
where \(N\) is already chosen and fixed. Choose
\begin{equation}%
\label{ChEps}%
\varepsilon=\kappa/3\,.
\end{equation}%
Let \(A(\kappa)=A(\kappa/3,N(\kappa))\), where \(A(\varepsilon,N)\) is the
value which appears in the formulation of the above stated Theorem which is attributed
to the names Cloizeaux-Mehta and Landau-Widom. Then for any \(a\) satisfying the
inequality \(a\geq{}A(\kappa)\), there exists the eigenvalue \(\sigma_{k(a,b(I))}\)
of the operator \(\mathscr{F}_{[-a,a]}^\ast\mathscr{F}_{[-a,a]}\) which belongs
to the interval \(I\). What is important that \emph{the value \(A(\kappa)\) depends only of
the length \(\kappa\) of the interval \(I\) but not on the position of \(I\) within the interval \([0,1]\).}

Thus, the following result is proved:
\begin{lemma}
\label{MaLe}%
Given a positive number \(\kappa\), there exists the number~\(\,\,A(\kappa)\),\,
\(A(\kappa)<\infty\),
 such that for every \(a\) satisfying the inequality \ \ \(a\geq~A(\kappa)\) the set
 of eigenvalues \(\big\lbrace\sigma_k(a)\big\rbrace_{k=0,\,1,\,2,\,\ldots\,\,}\) of the operator \(\mathscr{F}_{[-a,a]}^\ast\mathscr{F}_{[-a,a]}\)
 forms a \(\kappa\)-net for the interval \([0,1]\).
\end{lemma}

To formulate the result which is related to the operator \(\mathscr{F}_{[-a,a]}\)
itself rather to the operator \(\mathscr{F}_{[-a,a]}^\ast\mathscr{F}_{[-a,a]}\)
we  need a little bit modify the above reasoning.
\begin{lemma}
\label{MaLeM}%
Given a positive number \(\kappa\), there exists the number~\(\,\,A(\kappa)\),\,
\(A(\kappa)<\infty\),
 such that for every \(a\) satisfying the inequality \ \ \(a\geq~A(\kappa)\)
 each of the four sets
 of eigenvalues \(\big\lbrace\sigma_{4l+r}(a)\big\rbrace_{l=0,\,1,\,2,\,\ldots\,\,}\)
 of the operator \(\mathscr{F}_{[-a,a]}^\ast\mathscr{F}_{[-a,a]}\),
 corresponding to the residues
 \(r=0,\,1,\,2,\,3\), form a \(\kappa\)-net of the interval \([0,1]\).
\end{lemma}
\begin{proof}
For given \(a\), we can not control the residue class by\(\mod 4\) of the
value \(k(a,b(I)\). However we can control this residue class perturbing
a little bit the
value \(b(I)\). Given \(a,\,a\geq{}A(\kappa)\), and given one of the residues \(r,\,r=0,\,1,\,2,\,3\), we need to choose the value \(b\) such that
the following condition are satisfied:
\begin{subequations}
\label{Mod}
\begin{gather}
k(a,b)=r\ \ \mod4\,
\label{Mod1}\\
|(1+e^{b})^{-1}-(1+e^{b(I)})^{-1}|<\kappa/6\,,
\label{Mod2}
\\
-N\leq{}b\leq{}N\,.
\label{Mod3}
\end{gather}
\end{subequations}
Starting from \eqref{IntAs1}, we choose \(b\) such that the value
\((b-b(I))\ln{}a\) takes one of the seven values \(-3,\,-2,\,-1,\,0,\,1,\,2,\,3\).
Then \(|b-b(I)|\leq3/\ln{}a\). Thus, increasing \(A(\kappa\) if this is needed,
we can ensure the inequality \eqref{Mod2}.
Choosing \(b\) properly, we can ensure the equality \eqref{Mod1}.
Moreover, choosing the sign of the difference \(b-b(I)\) properly,
we can preserve the inequality \eqref{Mod3}.
(Actually, the last step is an overcautiousness.)
The inequality \eqref{Mod2} together with the inequality
\eqref{appr} ensure the inequality
\begin{equation}
\label{apprm}
\big|\sigma_{k(a,b)}-(1+e^{\pi{}b})^{-1}\big|\leq\kappa/2\,,
\end{equation}
(Recall that \(\varepsilon<\kappa/3\)\,.)

Thus if \(a\geq{}A(\kappa\), and \(r=\) takes one of the values \(0,\,1,\,2,\,3\), then in every interval \(I\) of length
\(\kappa\), \(I\subset[0,1]\), there exists an eigenvalue \(\sigma_k(a)\)
of the operator \(\mathscr{F}_{[-a,a]}^\ast\mathscr{F}_{[-a,a]}\) which
belongs to this interval: \(\sigma_k(a)\in{}I\), and the congruence \(k=r\,\,(\!\!\!\!\mod4)\)
holds.
\end{proof}
An immediate consequence of Lemma \ref{MaLeM} is the following
\begin{theorem}
\label{apprcr}
Given a positive number \(\kappa\), there exists the number \(\,\,A(\kappa)\),\,
\(A(\kappa)<\infty\),
 such that for every \(a\) satisfying the inequality \ \ \(a\geq{}A(\kappa)\)
 the set \(\big\lbrace{}\lambda_k(a)\big\rbrace_{k=0,\,1,\,2,\,}\) of the  eigenvalues
 of the truncated Fourier operator \(\mathscr{F}_{[-a,a]}\) forms an \(\kappa\)-net
 in the cross with the vertices \(\lambda=1\), \(\lambda=i\), \(\lambda=-1\), \(\lambda=-i\)\,.
\end{theorem}
\begin{proof}
If \(\sigma_k(a)\) is an eigenvalue of the operator
\(\mathscr{F}_{[-a,a]}^\ast\mathscr{F}_{[-a,a]}\), then the
number \(\lambda_k(a)=i^k(\sigma_k(a))^{1/2}\) is an eigenvalue of the operator
\(\mathscr{F}_{[-a,a]}\). If \(a\) chosen such that the numbers
\(\big\lbrace{}\sigma_{4l+r}(a)\big\rbrace_{l=0,\,1,\,2,\,}\)
form a \(\kappa^2\)-net of the interval \([0,1]\),  then the
numbers \(\big\lbrace{}\lambda_{4l+r}(a)\big\rbrace_{l=0,\,1,\,2,\,}\)
form a \(\kappa\)-net in the interval \([0,i^r]\).
\end{proof}
\begin{remark}
\label{NoCont}
Theorem \ref{apprcr} claims that for \(a\to\infty\), the set of eigenvalues
\(\big\lbrace{}\lambda_{k}(a)\big\rbrace_{k=0,\,1,\,2,\,}\) of the operator
\(\mathscr{F}_{[-a,a]}\)
fills in a sense the cross with the vertices \(\lambda=1\), \(\lambda=i\), \(\lambda=-1\), \(\lambda=-i\)\,. However, the spectrum of the limiting operator
\(\mathscr{F}_{[-\infty,\infty]}\), which is the non-truncated Fourier operator,
consists of the endpoints \(\lambda=1\), \(\lambda=i\), \(\lambda=-1\), \(\lambda=-i\)
of the cross only.
\end{remark}

\vspace*{5.0ex}
\noindent
\begin{minipage}[h]{0.45\linewidth}
Victor Katsnelson\\[0.2ex]
Department of Mathematics\\
The Weizmann Institute\\
Rehovot, 76100, Israel\\[0.1ex]
e-mail:\\
{\small\texttt{victor.katsnelson@weizmann.ac.il}}
\end{minipage}
\vspace*{3.0ex}

\noindent
\begin{minipage}[h]{0.45\linewidth}
Ronny Machluf\\[0.2ex]
Department of Mathematics\\
The Weizmann Institute\\
Rehovot, 76100, Israel\\[0.1ex]
e-mail:\\
\texttt{ronny-haim.machluf@weizmann.ac.il}
\end{minipage}
\end{document}